\documentclass[a4paper,oneside,12pt,reqno]{amsart}

\usepackage{amsmath}
\usepackage{amssymb}
\usepackage{latexsym}

\newtheorem{thm}{Theorem}
\newtheorem{defn}[thm]{Definition}

\newtheorem{cor}[thm]{Corollary}
\newtheorem{conj}[thm]{Conjecture}
\newtheorem*{theorem*}{Theorem}
\newtheorem{lem}[thm]{Lemma}
\makeatletter \let\@eqnsel = \hfil \makeatother 

\begin{document}

\title[On Abelian 3-d Difference Sets]{On Abelian Difference Sets  with 
  Parameters of 3-dimensional Projective Geometries}
\author{Kevin Jennings \\
UCD School of Mathematical Sciences\\
University College, Dublin}

\email{kevin.jennings@ucd.ie}


\begin{abstract}
A  difference set is said to have classical parameters if \linebreak
\[
(v,k, \lambda ) = \left(\frac{q^d-1}{q-1} , \frac{q^{d-1}-1}{q-1} ,
\frac{q^{d-2}-1}{q-1} \right).
\]
 The case $d=3$ corresponds to
planar difference sets.
We focus here
on the family of abelian difference sets with $d=4$. 
The only known examples of such difference sets 
correspond to the projective geometries $PG(3,q)$.
We 
consider an arbitrary difference set with the parameters 
of $PG(3,q)$ in an abelian group and
establish
constraints on its structure. In particular, we discern embedded substructures.
\begin{theorem*} 
Let $D$ be a normalized difference set with classical parameters in an
abelian group $G$ of order $(q^s+1)(q^{2s}+1)$, where $s$ is an odd
prime with $s \ge q$ and where $s \nmid q^2+1$.
Let $M$ be a subgroup of $G$ of order $(q+1)(q^2+1)$. 
Then $D \cap M$ is a normalized
difference set
with classical parameters
 in $M$.
\end{theorem*}

\end{abstract}

\maketitle
\section{Group Rings}
Let $G$ be a finite abelian group of order $v$ and let $\mathbb{Z}G$
denote the integral group ring of $G$.
Given an element $a = \sum a_ig_i \in \mathbb{Z}G$, we set
\[
a^{(-1)}  = \sum a_i g_i^{-1}.
\]
Let $D$ be a $k$-subset of $G$, where $k\ge 1$.
By a standard abuse of notation we will use the letter $D$ to
represent both the set of elements $D$ and the corresponding group ring element $D =
\sum_{d \in D} d$.

We say that $D$ is a 
$(v,k, \lambda)$-difference set in $G$ if $D$ satisfies the group ring equation
\[
DD^{(-1)} =\lambda G +n1, 
\]
where $ n: = k - \lambda$ is the \emph{order} of the difference set.
If $\lambda = 1$ the difference set is called
\emph{planar} and is associated with a projective plane of order $n$. 

It is elementary to show that $v,k, \lambda$ are related by 
the fundamental equation
\[
\lambda (v-1) = k(k-1)
\]
or equivalently
\begin{equation} \label{fund}
\lambda v = k^2 -n.
\end{equation}


Let $G$ be a finite group and let $H$ be a subgroup  of index $r$ in $G$. 
Let $D$ be a $(v,k,\lambda)$-difference set in $G$.
For $i = 1,2, \ldots, r$ and distinct cosets $Hx_i$ of $H$ in $G$, we define 
the \emph{intersection numbers}, $s_i$, of $D$ with respect to $H$ by
\[
s_i: = |D \cap Hx_i|.
\]
The following equations hold:
\begin{eqnarray}
 \sum_{i=1}^r s_i = k \label{interno1} \\
\sum_{i=1}^r s_i^2  =  \lambda |H| +n. \label{interno2}
\end{eqnarray}

We begin by drawing attention to an elementary  result on the 
distribution of elements of $D$ through cosets of any subgroup of $G$.

\begin{thm} \label{distribution}
Let $H$ be a subgroup of $G$ with $[G:H] =r$. Let $D$ be a $(v,k,
\lambda)$-difference set in $G$. Suppose that in a certain coset of
$H$ there are $s$ elements of $D$. Then
\begin{equation} \label{si}
\left| s - \frac{k}{r} \right| \leq \sqrt{n}\left( \frac{r-1}{r} \right).
\end{equation}
\end{thm}

\begin{proof}
We can take $s_r = s$ in equations \eqref{interno1} and
\eqref{interno2}. Then these read
$
\sum_{i=1}^{r-1} s_i = k - s
$
and
$
\sum_{i=1}^{r-1} s_i^2 = \lambda |H| +n - s^2.
$
The Cauchy-Schwarz inequality tells us that
$
\left( \sum^{r-1}_{i=1} s_i \right)^2 \leq (r-1) \left( \sum_{i=1}^{r-1} s_i^2 \right)
$
and hence
\[
(k-s)^2 \leq (r-1)( \lambda |H| +n -s^2).
\]
Multiplying out and 
noting that $v = r|H|$ and using \eqref{fund}
this simplifies to
\[
rs^2 - 2ks +2n - rn +\lambda |H| \leq 0.
\]
Completing the square 
and using \eqref{fund}
again, we obtain
\[
(s - \frac{k}{r})^2 \leq \frac{n(r-1)^2}{r^2}
\]
and hence the result.
\end{proof}

We observe that the mean distribution of $k$ elements across the $r$
cosets is $\frac{k}{r}$ so we can see from \eqref{si} that the number of
elements of $D$ in any coset is within $\sqrt{n}$ of this mean value.

We say that an automorphism $\sigma$ of $G$
is a  \emph{multiplier} for $D$ if $\sigma(D) = gD$ for some $g \in
G$. More particularly, we say $\sigma$ is a (numerical) multiplier
for $D$ if $\sigma(x)=x^m$ for all $x\in G$, where $m$ is an integer
relatively prime to $v$.
We call $gD$ a \emph{translate} of the difference set $D$. Clearly
$gD$ is itself a difference set.
We say that the difference set $D$ is \emph{normalized} if 
\[
\prod_{d
  \in D} d =1.
\]
 It is elementary
to show that any  difference set with $\gcd(v,k)=1$ has a unique translate which is
normalized. It is straightforward to prove that such a normalized difference set is fixed set-wise by any
numerical multiplier of $D$.

We quote now a well known refinement of the multiplier theorem of Marshall Hall. The proof can be found in
\cite[Section VI.4]{bjl}.

\begin{thm}[M. Hall] \label{mult}
Let $D$ be an abelian $(v,k, \lambda)$-difference set where $n = k -
\lambda$ is a power of a prime $p$ and $\gcd(p,v)=1$. Then the mapping $\sigma: x \to
x^p$ is a  multiplier for $D$.
\end{thm}

We also state an abridged version of the Mann Test. A proof of the general test can be found in \cite{bjl}.
\begin{thm}[Mann Test] \label{mann}
Let $D$ be a $(v,k, \lambda)$-difference set in an abelian group $G$ of order
$v$. Let $U$ be a subgroup of $G$ and let $G/U$ have exponent
$u^*$. Suppose $p$ is a prime not dividing $u^*$ and $p^f \equiv -1
\mod u^*$ for some $f \in \mathbb{N}$. Then the following hold:

\begin{enumerate}
\item   $n=p^{2j}n'$, where
$\gcd(p,n')=1$, for some $j \in \mathbb{Z}$.

\item For all
cosets $Ug$ of $U$ in $G$, the corresponding
intersection numbers $|D\cap Ug|$ of
$D$ relative to $U$ are congruent modulo $ p^j$.

\item $p^{j} \leq |U|$.
\end{enumerate}
\end{thm}

\section{Difference Sets with Classical Parameters}

Suppose now that $D$ is a $(v,k,\lambda)$-difference set  
with parameters \linebreak
\begin{equation} \label{classic}
(v,k, \lambda) = \left( \frac{q^d-1}{q-1}, \frac{q^{d-1}-1}{q-1},
  \frac{q^{d-2}-1}{q-1} \right),
\end{equation}
where $q$ is a power of a prime.
Any difference set with these parameters is said to have
\emph{classical parameters}.
The order of these difference sets is $q^{d-2}$ so
Theorem \ref{mult} applies.

If $d=3$ the difference set is said to be planar and this is the case
which has received the most interest. 
For the rest of this paper, $G$ will be an abelian group supporting a difference set
with classical parameters \eqref{classic} with $d=4$. So
\begin{equation} \label{d=4}
|G| = \frac{q^{4s}-1}{q^s-1} = (q^s+1)(q^{2s}+1)
\end{equation}
where $q^s$ is a power of a prime number.

We let $H$ be a subgroup of $G$ of order $q^s+1$. 
It is known in the folklore of the subject that $D$ has $2$-valued
intersection sizes with the cosets of $H$.
We gather here some useful results
 on abelian difference sets with classical parameters which apply to the
 case under investigation. The proofs rely on an application of the 
 Mann Test and can be found  in \cite[Theorems 3,4,7]{kpj}

\begin{thm} \label{kj}
Let $D$ be a normalized difference set with classical parameters in an
abelian group $G$ where $|G| = \dfrac{q^{4s}-1}{q^s-1}$ and let $H$ be
a subgroup of $G$ of order $q^s+1$, so $[G:H] = q^{2s}+1$. Then the
following hold:
\begin{enumerate}
\item $H$ is the unique subgroup of $G$ of order $q^s+1$;

\item $Syl_2(G)$ is cyclic, with generator $z$, say;

\item If $q$ is odd, then $Hz \subseteq D$ and $|D \cap Hx| =1$ for
  any other coset $Hx$ of $H$;

\item If $q$ is even, then $H \subseteq D$ and $|D \cap Hx| =1$ for
  any other coset $Hx$ of $H$. 
\end{enumerate}
\end{thm}

We note that by \eqref{si},
\[
|D\cap H| \leq \frac{k}{r} + \sqrt{n}\left( \frac{r-1}{r} \right)
\]
for any subgroup $H$ of $G$ of index $r$ where $G$ contains a
$(v,k,\lambda)$-difference set. 
In our case here, we have equality since a simple check verifies that
\[
\frac{k+\sqrt{n}(r-1)}{r} =
\frac{q^{2s}+q^s+1+q^s(q^{2s}+1-1)}{q^{2s}+1}
= q^s+1 = |H|.
\]
It is an interesting observation that $H$ is the largest  subgroup that can possibly be
contained inside a non-trivial difference set. The
spread of the elements of $D$ through the cosets of $H$ is as close to
the mean value $\frac{k}{r}$ as possible, so the distribution is unbiased.


\section{Singer Difference Sets and PG$(3,q)$}

Let $K = \mathbb{F}_q$ be the finite field with $q$ elements and let
$F= \mathbb{F}_{q^d}$. We view $F$ as a vector space over $K$ and
 let $K^*$ denote the multiplicative group of
 nonzero elements of $K$.
  Let $\pi$ be that natural epimorphism
  $\pi:F^* \to F^*/K^*$.
  Let $H$ be a
$K$-hyperplane in $F$. 
Then it is well known that $D:= \pi(H \backslash \{0\})$
is a Singer difference set
in $F^*/K^*$. 

We often choose $H$ to be the elements of
trace zero in $F$.
If $M$ is an intermediate field between $K$ and $F$,
we let $Tr_{F/M}$ denote the trace form from $F$ into $M$. 
$Tr_{M/K}$ is similarly defined.
If $\sigma$ denotes a power of the
Frobenius mapping $\sigma: x \to x^q$, then the
Galois group of $M$ is generated by an appropriate power of $\sigma$.
The trace of any element is given by the sum of its Galois
conjugates. 

\begin{thm} \label{extn} 
Let $K = \mathbb{F}_q$ be the finite field with $q$ elements. Let $F$
be the extension field of $K$ of degree $ab$, where $a,b \in
\mathbb{N}$. Let  $M$ be the
extension field of degree $a$ and $N$ be the extension field of degree $b$
over $K$. 
Let $D$ be the $M$-hyperplane in $F$ defined
by 
\[
D = \{ x \in F: Tr_{F/M}(x) = 0 \}
\]
and let $E$ be the $K$-hyperplane in $N$ given by 
\[
E = \{x \in N: Tr_{N/K}(x) = 0 \}.
\]
If $\gcd(a,b)=1$ then $E \subseteq D$ as subsets of $F$.
\end{thm}
\begin{proof}
It is enough to determine when we have
$Tr_{F/M}(x) = Tr_{N/K}(x)$, for an element $x \in N$.
In our situation, 
\[
Tr_{F/M}(x) = x + x^{q^a} +x^{q^{2a}} + \ldots + x^{q^{(b-1)a}},
\] 
 while
\[
Tr_{N/K}(x) = x+ x^q + x^{q^2} + \ldots + x^{q^{b-1}}.
\]
Let $x \in N$. Then $x^{q^b} = x$. Write 
$a = sb+t$,
where $s,t \in
\mathbb{Z}$ and  $ 0 \leq t \leq b-1$.
Then 
\[
x^{q^a} = x^{q^{sb+t}} = x^{q^{sb}q^t} = (x^{q^{sb}})^{q^t} = x^{q^t}.
\]
Considering $x$ as an element of $F$, 
\begin{eqnarray*}
Tr_{F/M}(x) & = & x+ x^{q^a} + \ldots + x^{q^{a(b-1)}}  \\
            & = & x+ x^{q^t} + \ldots + x^{q^{t(b-1)}} 
\end{eqnarray*}
where $t, 2t, \ldots , (b-1)t$ are all interpreted modulo $b$.
If the residues $t, 2t, \ldots , (b-1)t$ are all distinct modulo $b$, then the above
expression becomes
\begin{eqnarray*}
& = & x +x^q + \dots + x^{q^{b-1}}  \\
& = &  Tr_{N/K}(x). 
\end{eqnarray*}
If $\gcd(a,b) = 1$ then these residues are distinct and the result follows. 
\end{proof}

In particular,
if $a=4$ and $b$ is odd the condition in Theorem \ref{extn} is satisfied and 
$Tr_{F/M}(x) = Tr_{N/K}(x)$ for any $x \in N$. 
This gives us the following:

\begin{cor} \label{singer}
Let $D$ be the Singer difference set in a group $G$ with parameters 
$\left( \frac{q^{4s}-1}{q^s-1}, \frac{q^{3s}-1}{q^s-1}, \frac{q^{2s}-1}{q^s-1} \right)$
 derived from the trace zero hyperplane, where $s$ is odd.
Let $R$ be the subgroup of $G$ of order $\frac{q^4-1}{q-1}$.
Then $D \cap R$ is a 
$\left( \frac{q^{4}-1}{q-1}, \frac{q^{3}-1}{q-1}, \frac{q^{2}-1}{q-1} \right)$-
Singer difference set in $R$.
 \end{cor}

\section{Abelian Difference Sets with Parameters of PG$(3,q)$}

We mention here that the Singer difference sets are the only known examples of 
difference sets with these parameters in abelian groups. 
It has not been proved that this is the only structure possible.
Certainly it is a difficult open question whether all abelian planar difference sets are 
equivalent to the planar Singer difference sets.
We will now
consider an arbitrary difference set with the parameters 
of $PG(3,q)$ in an abelian group and
establish
constraints on its structure.
In particular, we try and generalize Corollary \ref{singer} to
any abelian difference set with these parameters. 
Generalizing from cyclic groups to abelian groups requires some
slightly clumsy looking technical conditions in
our hypotheses.

\begin{lem} \label{Mfix}
Let $G$ be an abelian group with $|G|=(q^s+1)(q^{2s}+1)$, where $s$ is
an odd integer. Suppose that $G$ contains a difference set with classical
parameters.
 Let $\tau$ denote the (multiplier) automorphism $\tau: x \to
x^{q^4}$. Let $M$ be any subgroup of order $ (q+1)(q^2+1)$ in $G$.
Then $M \le G^{\tau}$, where $G^{\tau}$ denotes the subgroup of fixed
points of $\tau$. Furthermore, if $Syl_r(G)$ is cyclic for all
prime numbers $r$ dividing $b$ or $c$ where $b = \gcd (q+1,s)$ 
and $c = \gcd(q^2+1,s)$ then $M = G^{\tau}$.  
\end{lem}

\begin{proof}
We observe that as $s$ is odd we have
\[
q^s+1 = (q+1)(q^{s-1}-q^{s-2} + \ldots -q+1)
\]
and correspondingly
\[
q^{2s}+1 = (q^2+1)(q^{2(s-1)}-q^{2(s-2)} + \ldots -q^2+1).
\]
Now $x \in G^{\tau}$ if and only if $x^{q^4-1} =1$ \emph{i.e.}
$
x^{(q-1)(q+1)(q^2+1)} =1.
$
From this it is clear that $M \leq G^{\tau}$, because
$m^{(q+1)(q^2+1)} =1$ for all $m \in M$.
We also note that $\gcd(q^4-1, |G|) = (q+1)(q^2+1)$.

Since $G$ contains a  difference set with classical parameters,
 the Sylow $2$-subgroup
of $G$ is cyclic by Theorem \ref{kj}. 
Let $r$ be a prime
dividing $\gcd(q+1, q^{s-1}-q^{s-2}+ \ldots -q+1)$. Now, since  $r$
divides $q+1$, we deduce that 
$
q \equiv -1 \mod r,
$
and, since $r$ divides $q^{s-1}-q^{s-2}+ \ldots -q+1$, we see that
\[
(-1)^{s-1}-(-1)^{s-2} + \ldots -(-1)+1 \equiv 0 \mod r,
\]
from which we conclude that
$
s \equiv 0 \mod r.
$
Hence $r$ divides $s$. We observe similarly that if $r$ is a prime
dividing $\gcd(q^2+1, \dfrac{q^{2s}+1}{q^2+1})$, then again, $r$ divides $s$.

We let $b = \gcd (q+1,s)$ and $c = \gcd(q^2+1,s)$. 
If, for
all prime divisors $r$ of $b$ or of $c$,
each Sylow $r$-subgroup is cyclic
 then $M =
G^{\tau}$. 
We note that if $s$ is a large enough prime  
($s >q^2$ for instance), then $b = c =1$
and this
condition will be automatically satisfied.
Of course, if $G$ is cyclic then $G^{\tau} = M$.
\end{proof}

We observe  that the same power of $2$ divides both $|M|$ and $|G|$,
since 
$[G:M]$ is odd,
and hence $M$ contains the Sylow
2-subgroup of $G$, which is cyclic by Theorem \ref{kj}, generated by $z$.

As before, let $H$ be the subgroup of $G$ of order $|H| = q^s+1$.
We have that $Hz \subseteq D$ by Theorem \ref{kj}.

To generalize Corollary \ref{singer} we would like to show that $D
\cap M$ is a difference set for $M$. We first show,
subject to
weak restrictions, that $D \cap M$ has the correct size to be a
difference set with classical parameters.

\begin{lem} \label{size}
Let $D$ be a normalized difference set in an abelian group $G$ of
order $(q^s+1)(q^{2s}+1)$, where $s$ is odd. Let $M$ be a subgroup of $G$ with $|M| =
  (q+1)(q^2+1)$. Let $b = \gcd(q+1,s)$ and $c = \gcd(q^2+1,c)$. For each prime divisor $r$ of
  $b$ or of $c$, suppose that $Syl_r(G)$ is cyclic.
  Then $|D \cap M| = q^2+q+1$.
\end{lem}
\begin{proof}
The conditions on $Syl_r(G)$ being cyclic ensure that $M = G^{\tau}$ here, 
as in Lemma \ref{Mfix}.
Let $H$ be the subgroup of $G$ of order $q^s+1$.
As shown in the proof of Lemma \ref{Mfix}, the hypotheses ensure that
\[
\gcd(q+1, \dfrac{q^s+1}{q+1}) =1
\]
and hence we can express $H$ as a direct product
$
H = AB,
$
where $|A| = q+1$ and $|B| = \dfrac{q^s+1}{q+1}$
with $A \cap B = 1$.
By the above and the cyclicity of the Sylow $2$-subgroup, $A$
is the unique
subgroup of order $q+1$ in $G$.
Now, since $\gcd(|M|,|H|) = \gcd((q+1)(q^2+1),q^s+1) = q+1$, 
we can deduce that
$
H \cap M = A
$
and thus $H \cap M$ is the unique subgroup of order $q+1$ in $G$.
 
Now $|M:H \cap M| = q^2+1$. We can decompose $M$ as 
\begin{displaymath}
M = \bigcup^{q^2+1}_{i=1} (H \cap M)m_i  
\end{displaymath}
where $ m_1, \ldots ,  m_{q^2+1} $ are distinct coset representatives of
$H \cap M$ in $M$.
Then we claim that 
\[
\bigcup^{q^2+1}_{i=1} H m_i
\]
is a union of different cosets of $H$ in $G$. 
For  $Hm_i = Hm_j$ implies $m_im_j^{-1} \in H \cap M $, since it is
clearly in $M$, 
 and hence $ (H \cap M)m_i = (H \cap M)m_j$ from which $m_i = m_j$,
 proving the claim.

As observed above,
$M$ contains the Sylow $2$-subgroup of $G$. If $z$ is a generator for
$Syl_2(G)$  then 
$z \in M$ and
$ z \notin H \cap M$, since $H$ does not contain $Syl_2(G)$.
Choose $m_1 = z$. Since $Hm_1 \subset D$, we have that 
\[
(H \cap M) m_1 \subset D.
\]
Now each of the distinct cosets $Hm_i$ for $i \geq 2$  contains a unique element
of $D$, by Theorem \ref{kj} again.
So for each such representative $m_i$, there exists a unique $h_i \in H$
with
$
h_im_i \in D.
$
The multiplier $\tau: x \to x^{q^4}$ fixes $M$ by Lemma \ref{Mfix} and
\[
\tau(h_i) \tau(m_i) \in \tau(D) = D,
\]
and thus
\[
 \tau(h_i)m_i \in D.
\]
Also, $\tau(h_i) \in H$, so by the uniqueness of $h_i$ we must have
$
\tau(h_i)=h_i.
$
Since $M = G^{\tau}$ here,
we deduce that $h_i \in M$ for each $i$ and thus
$h_im_i \in D \cap M$ for $i \geq 2$.
We conclude that
\[
D \cap M = (H \cap M) m_1 \cup \{h_i m_i :i=2,3,\ldots , q^2+1  \}.
\]
Hence $|D \cap M| = (q+1) +  q^2 = q^2+q+1$.
\end{proof}



 The following theorem is our generalization of Corollary
\ref{singer} to abelian difference sets.
Apart from 
Jungnickel and Vedder's result on
 planar difference sets  with square order
(Theorem \ref{rod4} here),
this is the only case we know of 
 where the 
parameters of the difference set guarantee a subdifference set.

\begin{thm} \label{main2}
Let $D$ be a normalized difference set with classical parameters in an
abelian group $G$ of order $(q^s+1)(q^{2s}+1)$, where $s$ is an odd
prime with $s \ge q$ and where $s \nmid q^2+1$.
Let $M$ be a subgroup of $G$ of order $(q+1)(q^2+1)$. 
Then $D \cap M$ is a normalized
difference set
with classical parameters
 in $M$.
\end{thm}

\begin{proof}
Our hypothesis that 
$s \nmid q^2+1$ guarantees,
by Lemma \ref{Mfix},
that $M$ is the group of fixed points of the multiplier $\tau: x \to
x^{q^4}$. The orbits of $\tau$ have length dividing $s$
and since $s$ is a prime number,
the orbits of $\tau$ have length $1$ or $s$.
The orbits of length $1$ correspond to elements of $M$, since $M = G^{\tau}$.
Let $g \in M$. 
Then there exist $\lambda = q^s+1$ ordered pairs 
$(a_i,b_i) \in D \times D$ with $g = a_ib_i^{-1}$.

Now 
\[
\tau(g) = g = \tau(a_i) \tau(b_i^{-1}) = \tau(a_i) \tau(b_i)^{-1} 
\]
and
\[
(\tau(a_i),\tau(b_i)) \in D \times D.
\]
So the $\lambda$  ordered pairs representing $g$ come in multiplier orbits of
length $1$ or $s$.
Now, since  $s$ is prime,  we have
\begin{equation} \label{cong}
\lambda = q^s+1 \equiv q+1 \mod s 
\end{equation}
and we see that if $s>q+1$,
we must have $q+1$ ordered pairs which are fixed by $\tau$.
Thus each element of $M$ can be represented by exactly $q+1$
pairs $(a_i,b_i) \in D \cap M \times D \cap M$.

Indeed, even if $s=q$ is a prime or in the case where $s = q+1$ and $s$ is a
Mersenne prime, the conclusion remains valid. To show this, we observe that
each element of $M$ can be represented at least twice as a
``difference'' from $D \cap M \times D \cap M$. This is because 
\[
D \cap M = (H \cap M) m_1 \cup \{h_i m_i :i=2,3,\ldots , q^2+1  \}
\]
from the arguments in Lemma \ref{size} 
and the products $ab^{-1}$ where $a \in (H \cap M)m_1$ and $b \in \{\cup
h_i m_i \} $ cover each  element of $M$ as do the products
$ab^{-1}$ where $a \in \{\cup h_i m_i \} $ and  $b \in (H \cap M)m_1$.
So \eqref{cong} will guarantee the result if $s>q-1$.
\end{proof}

\section{Another Subgroup with two-valued intersection numbers}

In our case of difference sets with classical parameters
 where $d=4$,
 we have a second subgroup whose intersection numbers are two-valued.

\begin{thm} \label{DintK}
Let $G$ be an abelian group with 
$
|G| = (q+1)(q^2+1),
$
where $q$ is a power of a prime, and let $D$ be a normalized difference set with classical parameters in $G$. 
Then,
\begin{enumerate}

\item There is a unique subgroup $K$ of order $q^2+1$ in $G$.

\item
$
|D \cap Kx| = \begin{cases}
1 & \textit{for one distinguished coset} \\
q+1 &  \textit{for the other } q \textit{ cosets}.
\end{cases}
$ 
\end{enumerate}
\end{thm}

\begin{proof}
Firstly, if $q$ is even then the Sylow $2$-subgroup of $G$ is trivial, while if $q$ is odd then Theorem \ref{kj}
tells us that the Sylow $2$-subgroup is cyclic. Since 
$$\gcd(|K|,|G:K|) = \gcd(q^2+1,q+1) $$ 
is a divisor of $2$, we deduce that $K$ is unique.

We note that $|G:K|= q+1$ and that $K$ satisfies the role of $U$ in our statement of the Mann Test, 
Theorem  \ref{mann}.
Letting $s_i$ denote the intersection number $|D \cap Kx_i|$, we deduce from the Mann Test part (c) that 
all the $s_i$ are congruent to each other modulo $q$, say $s_i \equiv y \mod q$, where $0\le y<q$.
 As before, we have
\begin{equation} \label{k}
\sum^{q+1}_{i=1} s_i = k = q^2+q+1
\end{equation}
and hence
\[
(q+1)y +rq = q^2+q+1,
\]
for some $r \in \mathbb{Z}$. 
It is straightforward to see that $y = 1$ and hence $r=q$.
Since there are $q+1$ cosets,
we conclude that at least one coset of $K$ has intersection size $1$ with $D$ (since otherwise
all $s_i$ satisfy $s_i\ge q+1$, which is impossible from \eqref{k}).

Finally, the Cauchy-Schwarz inequality completes the proof.
Recall that 
\begin{equation} \label{CS}
\left( \sum_{i=1}^n a_i \right)^2 \leq n \left( \sum_{i=1}^n a_i^2 \right),
\end{equation}
with equality if and only if all the $a_i$ are equal.
Letting $s_{q+1}=1$, we have that
\[
\sum_{i=1}^q s_i = (q^2+q+1) -1 = q^2+q
\] 
while \eqref{interno2} yields
\[
\sum_{i=1}^{q} s_i^2 = ((q+1)(q^2+1)+q )- 1 = q(q+1)^2 .
\]
We now have equality in \eqref{CS}
and hence all the $s_i$ must be equal to each other.
By counting, $s_i = q+1$ for all the other intersection numbers.
\end{proof}

We note that since $|G:K| = q+1$,
the multiplier $\sigma:x \to x^q$ is an involution on the cosets of $K$.
Now the coset $Kx$ is fixed by $\sigma$ if $Kx = Kx^q$
which occurs if and only if $x^{q-1} \in K$. 
Since $\gcd(q-1,q+1)$ divides $2$ and since the Sylow $2$-subgroup
of $G$ is cyclic, the only cosets fixed by $\sigma$ are $K$ and $Kw$,
where $w^2 \in K$ but $w \notin K$. 
Since $|K| \equiv 2 \mod 4$, the unique element of order $2$ in $G$ is in $K$.
Hence $w$ must be an element of order $4$ in $G$
and the distinguished coset is either $K$ itself or $Kw$.

In particular, when $q$ is even, there is only one coset, $K$ itself,
fixed by $\sigma$. Hence $D \cap K = \{1 \}$ in this case. We summarize 
this in the following:

\begin{cor} \label{HK}
Let $D$ be a normalized $(v,k, \lambda)$-abelian difference set in $G$ with parameters \eqref{d=4}. 
Suppose $q$ is even, say $q=2^s$ so that $|G|= (2^s+1)(2^{2s}+1) $ and $G$ is a direct
product $G=HK$ where $|H| = 2^s+1$ and $|K| = 2^{2s} +1$. Then $H \subseteq D$ and
  $|D \cap Hx|=1$ for each other coset of $H$.
Furthermore $D \cap K = \{ 1\}$ and $|D \cap Kx| = 2^s+1$ for each other coset of $K$.
\end{cor}


\section{Minimal Difference Sets and Conjectures}

\begin{thm} \label{1573}
Let $D$ be a normalized difference set in an abelian group $G$ with parameters
\[
\left( \frac{2^{4s}-1}{2^s-1}, \frac{2^{3s}-1}{2^s-1} , \frac{ 2^{2s}-1}{2^s-1} \right).
\]
Suppose $s$ is odd. Then  $G$ has a subgroup, $M$, of order $15$ and
 $D \cap M$ is a $(15,7,3)$-difference set in $M$.
\end{thm}
\begin{proof}
As in Corollary \ref{HK}, we write 
$
G = HK,
$
where $|H| = 2^s+1$ and $|K| = 2^{2s}+1$.
Since $s$ is odd, we have that
$
2^s+1 \equiv 0 \mod 3 
$ 
and
$
2^{2s}+1 \equiv 0 \mod 5.
$
So $3$ divides $|H|$ and $5$ divides $|K|$.
Let $k \in K$ have order $5$ in $K$. Then since $|D \cap Hk|=1$, by Corollary \ref{HK},
there exists a unique $h \in H$ with $hk \in D$.

By Theorem \ref{mult}, $\sigma: x \to x^2$ is a multiplier fixing
$D$.
So 
\[
\sigma^4(hk)  =  \sigma^4(h)\sigma^4(k)  =  \sigma^4(h) k \in  D
\]
since $\sigma$ fixes $D$. Now, since $h$ is unique,
$ \sigma^4 (h)  =  h$
and thus
$ h^3  =  1$.
Finally $h \neq 1$ since $D \cap K = 1$, by Corollary \ref{HK}.

Let $M=<hk>$ be a subgroup of $G$ of order $15$. Then 
\[
D \cap M = \{ 1, h, h^2, hk,h^2k^2,hk^4,h^2k^3 \}
\]
which is a $(15,7,3)$-difference set for $M$. It can be seen directly
in this case that 
 each element of $M$ arises exactly $3$ times as a difference from this set.
\end{proof}

We have called this $(15,7,3)$-difference set a \emph{minimal difference set}
as a copy of it appears embedded in the structure of larger members 
of the family.
We feel that it recalls the role of the prime subfield in field theory and
that it also echoes Ho's result in Theorem \ref{hothm},
where a Baer subplane is embedded in the 
structure of larger members of the family .

Amongst the results on planar abelian difference sets which 
motivated this work, we highlight
the following theorems: the first due  to Ostrom \cite{ostrom}
in the cyclic case and then extended to the abelian case by Jungnickel and Vedder \cite{jvedder};
and the second due to Ho \cite{ho}.

\begin{thm}[Jungnickel and Vedder] \label{rod4}
 Let $G$ be a finite abelian
group and let $D$ be a normalized planar difference set
of square order $m^2$ in $G$. Let $H$ be the
unique subgroup 
of order $m^2+m+1$ in $G$. Then $D\cap H$ is a 
normalized planar difference set
of order $m$ in $H$.  
\end{thm}

\begin{thm}[Ho] \label{hothm}
Let $D$ be a planar difference set of order $m^s$ in the cyclic group $G$. 
Then $D$ contains a planar difference set of order $m$ in the 
unique subgroup of order $m^2+m+1$ of $G$ 
if and only if $s$ is not a multiple of $3$.
\end{thm}

\begin{defn}
We call a difference set $D'$ which has the parameters 
$((q+1)(q^2+1),q^2+q+1,q+1)$ a \emph{minimal difference set}
if $q=p^r$ where $p$ is a prime and $r$ is a power of $2$.
\end{defn}

We have a partial generalization of Corollary \ref{singer}  for
abelian difference sets in Theorem \ref{main2} but we are still a long way 
from the following conjecture:

\begin{conj}
Let $D$ be a normalized difference set with parameters \\
$((q^s+1)(q^{2s}+1),q^{2s}+q^s+1,q^s+1)$
in an abelian group $G$. Then $D$ contains a minimal difference set embedded in it,
in the sense that there exists a subgroup $S$ of $G$ with
$D\cap S =D'$, where $D'$ is a minimal difference set.
\end{conj}

It would be sufficient to prove this result  true for any odd prime $s$, 
in which case, a cascade effect would guarantee the result for any odd $s$. In
Theorem \ref{main2} we have shown that for given $q$, the conjecture
 is true for all  primes $s$
larger than $q^2$. 
By Theorem \ref{1573}, it is true for all $s$ when $q=2$.


\begin{acknowledgements}
This work is part of the author's Ph.D. thesis. The author is very
grateful to Rod Gow (UCD) for his guidance and advice.
\end{acknowledgements}


\bibliographystyle{amsplain}
	\bibliography{bibliography}

\providecommand{\bysame}{\leavevmode\hbox to3em{\hrulefill}\thinspace}
\providecommand{\MR}{\relax\ifhmode\unskip\space\fi MR }
\providecommand{\MRhref}[2]{%
  \href{http://www.ams.org/mathscinet-getitem?mr=#1}{#2}
}
\providecommand{\href}[2]{#2}
\begin{thebibliography}{1}

\bibitem{bjl}
Thomas Beth, Dieter Jungnickel, and Hanfried Lenz, \emph{Design theory. {V}ol.
  {I}}, second ed., Encyclopedia of Mathematics and its Applications, vol.~69,
  Cambridge University Press, Cambridge, 1999. \MR{MR1729456 (2000h:05019)}

\bibitem{ho}
Chat~Yin Ho, \emph{Subplanes of a tactical decomposition and {S}inger groups of
  a projective plane}, Geom. Dedicata \textbf{53} (1994), no.~3, 307--326.
  \MR{MR1311323 (96i:51009)}

\bibitem{kpj}
K.~P. Jennings, \emph{On substructures of abelian difference sets with
  classical parameters}, Journal of Combinatorial Designs (to appear 2007).

\bibitem{jvedder}
Dieter Jungnickel and Klaus Vedder, \emph{On the geometry of planar difference
  sets}, European J. Combin. \textbf{5} (1984), no.~2, 143--148. \MR{MR753004
  (85i:05047)}

\bibitem{ostrom}
T.~G. Ostrom, \emph{Concerning difference sets}, Canadian J. Math. \textbf{5}
  (1953), 421--424. \MR{MR0056006 (15,10c)}

\end{thebibliography}
\end{document}